\documentclass[12pt]{amsart}
\usepackage{epsfig,amssymb,epic,eepic,color}
\usepackage{amsmath,amsthm,amsfonts,amssymb}
\newtheorem{thm}{Theorem}[section]

\newtheorem{prop}[thm]{Proposition}
\newtheorem{defn}[thm]{Definition}
\newtheorem{lemme}[thm]{Lemma}

\newcommand{\be}{\begin{enumerate}}
\newcommand{\ee}{\end{enumerate}}
\newcommand{\bi}{\begin{itemize}}
\newcommand{\ei}{\end{itemize}}

\def\R{\mathbb{R}}

\def\ga{\gamma}

\def\al{\alpha}
\def\be{\beta}

\def\vp{\varphi}

\def\la{\lambda}

\def\si{\sigma}

\def\ep{\varepsilon}

\def\ds{\displaystyle}

\def\nd{\noindent}
\def\bull{\hfill$\Box$}

\begin{document}
\thispagestyle{empty}
\vskip 1cm
\title{Positive Legendrian regular homotopies}
\author{Fran\c cois Laudenbach}
\address{Universit\'e de Nantes, UMR 6629 du CNRS,
44322 Nantes, France.}
\email{francois.laudenbach@univ-nantes.fr}

\keywords{Legendrian immersion, positive regular homotopy}
\subjclass{53D10}
\begin{abstract}  In contrast with what happens for Legendrian embeddings,
 there 
always    exist positive loops of Legendrian immersions.  
\end{abstract}
\maketitle
\vskip 1cm
\section{Introduction}

\medskip

In this note, the following is proved:

\begin{thm} \label{theo} Let $L$ be an n-dimensional closed manifold, 
$J^1(L,\R)$ be its
 space of 1-jets, $\al=dz-pdq$ be its canonical contact form and $\xi=\ker\al$
be its associated contact structure. There exists a loop of Legendrian 
immersions $\vp_t: L\to J^1(L,\R)$, $t\in S^1$, starting from the 0-section 
and positive in the following sense: $\ds \frac{\partial\vp}{\partial t}$ is
 transverse to $\xi$ at every point. Moreover, $\vp_t$ may be chosen
 $C^0$-close to the stationary loop.\footnote{As a consequence, the theorem
 holds for any contact manifold $(M,\xi)$ and any Legendrian immersed 
submanifold $L$ (apply the {\it tubular  neighborhood theorem} for Legendrian immersions in a contact manifold). I am grateful to the referee for this remark and, more generally, for  
his very careful reading of the first version.}\\
\end{thm}

In \cite{polt} Y. Eliashberg and L. Polterovich emphasized that the existence
 of positive contractible loops of Hamiltonian contact diffeomorphisms on a 
contact manifold should have an important topological significance.  
In \cite{ekp} with S. S. Kim, they
proved 
that such loops exist on the standard spheres $S^{2n+1},\ n>0,$ but not on 
$J^1(L,\R)$.
 By using invariants defined by C. Viterbo in \cite{viterbo}, V. Colin,
 E. Ferrand and P. Pushkar proved that there do not exist any positive 
loops 
of Legendrian embeddings of $L$ into  $J^1(L,\R)$
starting from  the zero-section (see \cite{colin}). 
Thereafter, Emmanuel Ferrand
 asked me about the same question in replacing embeddings by immersions; 
I am grateful to him
for that question. I thought 
that methods which are developed in the marvellous book by Y. Eliashberg and 
N. Mishachev \cite{eliash} could apply. This is the case as  is explained 
below; actually we do not apply a known  $h$-principle, but we prove the 
$h$-principle for the problem at hand. I thank 
Vincent Colin for interesting discussions on related topological questions.

I am glad to have here the opportunity to express my deep admiration to 
Yasha.\\

\section{Field of transverse loops}

\medskip
 Let $X$ be the  Reeb vector  field on $J^1(L,\R)$ associated to $\al$, 
that is the unique contact vector field such that $\al(\xi)=1$.
 A {\it transverse loop}
will be an immersed closed curve which is transverse to $\xi$.\\

\begin{prop}\label{loops}There exists  a smooth family of transverse loops
 passing through every
 point in $L$, identified with the 0-section of $J^1(L,\R)$. 
More precisely, there exists $\vp:L\times S^1\to J^1(L,\R)$ 
such that:
\begin{itemize}
\item[{\rm 1)}] $\vp(x,0)=x$ for every $x\in L$;
\item [{\rm 2)}] $\ds\frac{\partial\vp}{\partial t}(x,t)= 
X\bigl(\vp(x,t)\bigr)$
 when $t$ is close to 0 in $S^1$ (here $t$ stands for the variable in $S^1$);
\item [{\rm 3)}] $\ds\frac{\partial\vp}{\partial t}(x,t)$ is transverse 
to $\xi$ for every  $(x,t)\in L\times S^1$.
\end{itemize}
\end{prop}

\nd {\bf Proof.} Here we think of $S^1$ as $[-1,+1]$ with identified 
end points. If $X^t$ denotes  the flow of $X$, we define the required family by 
$\vp(x,t)=X^t(x)$ when $t\in [-1/2,+1/2]$. For  extending $\vp$ to 
$L\times S^1$, we will use a suitable family of contact Hamiltonians.

We  choose a  triangulation $\mathcal T$ of $L$ and a finite  open covering
 of $L$ by charts $\{\mathcal O_i\}_{i\in I}$ such that each simplex $\si$
is contained in some $\mathcal O_{i(\si)}$ and, when  $\tau $ is a face of $\si$, 
$\mathcal O_{i(\tau)}\subset\mathcal O_{i(\si)}$.
 Let $U_i=J^1(\mathcal O_i,\R)$. We are going
 to construct $\vp$ so that $\vp(\si\times S^1)\subset U_{i(\si)}$ for each 
$\si\in\mathcal T$.

We first describe a transverse loop passing through a vertex $A$. Let $A^\pm
=X^{\pm 3/4}(A)$. We choose a smooth 
path $\ga$ from $A^+$ to $A^-$ in $U_{i(A)}$  
avoiding $X^t(A),\ t\in [-1/2,+1/2]$; it is not required to be transverse 
to $\xi$. It may be viewed as the trajectory of 
$A^+$ during a contact isotopy $\psi_s,\ s\in [0,1]$,
generated by some time dependent contact 
Hamiltonian $h_s$. We assume that $h_s$  
is constant  on a small interval of  the Reeb orbit 
$[X^{-\ep}\bigl(\ga(s)\bigr), X^{+\ep}\bigl(\ga(s)\bigr)]$. The Hamiltonian is 
 extended to $J^1(L,\R)$ with  compact support in  $U_{i(A)}$ avoiding 
$X^t(A),\ t\in [-1/2,+1/2]$. The arc $\psi_1\bigl(X^t(A)\bigr), \ t\in 
[-1/2,+3/4]$ is transverse to $\xi$, starts at $X^{-1/2}(A)$ and 
 ends at $X^{- 3/4}(A)$. 
Moreover, it is vertical, {\it i.e.}
 tangent to $\ds\frac{\partial}{\partial z}$, near its end points.
Therefore this arc can be smoothly closed by adding the  vertical interval
$J(A)=\{X^t(A)\mid t\in [-3/4,-1/2]\}$. 
In general the resulting  loop is only
immersed with possible double points in $J(A)$.\\

Now, assume  recursively that 
$\vp$ is already constructed  on  $\partial\si\times S^1$, where 
$\si$ is a simplex of $\mathcal T$. We perform a similar construction where
$A$ is replaced by $x$ running in $\si$. We have to find a family of paths 
$\ga_x$ from $x^+$ to $x^-$ in  $ U_{i(\si)}$, avoiding 
$X^t(x),\ t\in [-1/2,+1/2]$, depending smoothly on $x$ and extending the 
given family on $\partial\si$. As the dimension of $ U_{i(\si)}$ is much bigger 
than the dimension of $\si$, there is no obstruction for solving this homotopy problem.
\bull

\medskip
\nd {\bf Remark.} By rescaling it is easy to make the diameter of each loop
 in the family $\vp$ less than any positive $\ep$.\\

In the next step, we will thicken each tranverse loop $\vp_x,\ x\in L,$ from
the family $\vp$ yielded by proposition \ref{loops}, into an immersed
 tube equipped 
with a contact action of $S^1$, depending smoothly on $x$. 
Even forgetting that
 the 
action has to be contact, such a family of thickenings is not
 so obvious to get, 
except when $\xi$ is a trivial bundle. Of course one tranverse loop may be 
thickened into a tube with a $S^1$-action, but a problem appears when one
 tries to do it in a family. For that, we will take care of the concept of
 trivialization.
For $x\in L$, we use the same notation $\vp_x$ for denoting the immersed curve 
$\vp(x, S^1)$ or its parametrization $t\in S^1\mapsto \vp(x,t)$. We think of 
$\nu_x:=\xi|\vp_x$ as its normal bundle, which is a trivial bundle since the 
loop lies  in a chart. Classically, a trivialization of $\nu_x$
is  a framing. But, if $\xi$ is not a  trivial bundle, such a framing does not 
exist if it is required to depend continuously on $x$. We use another 
definition of trivialization.\\

\begin{defn} Let $A$ be a submanifold of $J^1(L,\R)$, with a base point $a$.
A trivialization of $\xi|A$ is a smooth field of linear isomorphisms
$\pi_z:\xi(z)\to\xi(a),\ z\in A$.\\
\end{defn}

For instance, if $A$ is a convex Darboux chart, we have a canonical
trivialization obtained by parallelism and projection along $\ds\frac{\partial}
{\partial z}$. Therefore, if $\si$ is a simplex in $\mathcal T$ and 
$x\in\si$, the Darboux chart $U_{i(\si)}$ induces a canonical trivialization
on $\nu_x$ using $\vp(x,0)=x$ as a base point of $\vp_x$.
 When $\tau$ is a face of $\si$, for $x\in\tau$ we have two 
trivializations of $\nu_x$: one when $\vp_x$ is viewed as a loop in 
$U_{i(\si)}$ 
 and the other one when it is viewed as a loop in $U_{i(\tau)}$. But, if the 
triangulation is fine enough and the diameter of the loops $\vp_x$ is small
 enough, both trivializations are close to each other. Arguing in this way, we
 have the following.\\

\begin{lemme}\label{linear} There exists  a smooth family of trivializations 
 $\pi_{x,t}:\xi(\vp(x,t))\to\xi(x),\ (x,t)\in L\times S^1$.
 This family may be viewed as 
a family of linear actions of $S^1$ on $\nu_x$, depending smoothly
 on $x\in L$.\\
\end{lemme}

\nd{\bf Proof.} The question is to find a global section of some bundle 
$\mathcal E$ over 
$L\times S^1$ whose fiber at $(x,t)$ is the space of linear isomorphims
$\pi_{x,t}:\xi\bigl(\vp(x,t)\bigr)\to \xi(x)$. Two such isomophisms differ
 by an automorphim of $\xi(x)$ depending on $t$. Therefore we introduce  
the group bundle  $\mathcal A$ over $L\times S^1$ whose fiber
$A_{x,t}$ over $(x,t)$ 
is the group of linear automorphisms of 
$\xi(x)$. Around the identity section in $\mathcal A$
 we choose $n$ disk sub-bundles
$V_1\subset V_2\ldots \subset V_n$ such that, if $f$ and $g$ are in 
$V_j\cap A_{x,t}$, then $f^{-1}$ is in $V_j$ and $f\circ g$ is in 
$V_{j+1}$. As above, we choose the triangulation $\mathcal T$ fine enough 
and the diameters of the loops $\vp_x$ small enough 
so that  the following condition is fulfilled:
\bi
\item[ ]for each simplex $\si\in \mathcal T$ and any face $\tau$, and for 
 any
$x\in \tau$, the canonical trivializations of $\nu_x$ associated to the 
Darboux charts $U_{i(\si)}$ and $U_{i(\tau)}$ differ by a section of $V_1$ 
over $\{x\}\times S^1$.\ei
Now the wanted section of $\mathcal E$  will be constructed  recursively 
over the skeleta $\mathcal T^{[k]}\times S^1$. Let $\si$ be a $k$-simplex 
of $\mathcal T$. Assume that we have a smooth family of trivializations 
$\pi_{x,t}$ of $\nu_x$ for 
$x\in \partial\si$ and assume that, for any $(k-1)$-face $\tau$ of $\si$, 
it differs
from the canonical trivialization associated to $U_{i(\tau)}$ by a section 
of $V_{k-1}$. Hence it differs from the canonical trivialization associated 
to $U_{i(\si)}$ by a section of 
$V_1\circ V_{k-1} \subset V_k$. As the fiber of $V_k$ is contractible, 
two sections of $V_k$ are homotopic and the family $\{\pi_{x,t}\}$ extends 
over $\si\times S^1$.
\bull\\

\begin{prop}\label{tube} There exists a family of closed immersed tubes 
$\theta_x: D_x\times S^1\to J^1(L,\R), \ x\in L$, where $D_x$ 
is a small $2n$-dimensional disk centered at $x$ and tangent to 
$\xi(x)$ at $x$, with the following properties:
\begin{itemize}
\item[1)] $\theta_x$ depends smoothly on $x\in L$;
\item[2)] $\theta_x(x,t)=\vp(x,t)$ for every $(x,t)\in L\times S^1$;
\item[3)] for $z\in D_x$, $\theta_x(z,0)=\theta_x(z,1)=z$ and, 
for $t$ close to $0=1$ in $S^1$, 
 $\theta_x(z,t)=X^t (z)$; 
\item[4)]  $S^1$ acts on the source of $\theta_x$
 by contact diffeomorphims with respect to 
the contact structure $\theta_x^*(\xi)$.
\end{itemize}
\end{prop}

 \nd{\bf Proof.} According to lemma \ref{linear}, we have a family of linear
  $S^1$-actions on the $2n$-disk bundle $\nu_x^r$ of radius $r$
 about the 0-section
 in $\nu_x$, $x\in L$; here an $S^1$-invariant metric is used.
By choosing an exponential map $exp\,:\,\xi\to J^1(L,\R)$ and taking $r$ small 
enough, $exp(\nu_x^r)$ is a family of immersed tubes
$\theta_x: D_x\times S^1\to J^1(L,\R), \ x\in L$, which meets  the 
wanted conditions but the contactness; here $D_x$ is the fiber of the tube
 at $x$. The contact form $\al$ induces on $D_x$ a Liouville form $\la_x$ for
 the symplectic structure induced by $d\al$ on $D_x$.
Hence we get  an $S^1$-invariant  contact form on the tube:
$\tilde\al_x=dt+\la_x$. So we have two contact forms on $ D_x\times S^1$,
$\theta_x^*(\al)$ and $\tilde\al_x$. The underlying plane fields 
are both transverse to the $S^1$-orbits and coincide along
the core of the tube $\{x\}\times S^1$. Gray's stability theorem \cite{gray}
 applies relatively to the core curve
and yields a conjugation of the germs of both contact structures
along the core. Carrying the $S^1$-action over this conjugation we fulfill
 condition 4) on a small tube around   $\{x\}\times S^1$. \bull\\

We thicken the disk $D_x$ into a $2n+1$-ball $B_x$ contained in the cylinder 
$C_x:=\cup_{ s\in [-\ep,+\ep]}X^s(D_x)$. A point $z\in C_x$ reads $z=X^s(y)$, 
$y\in D_x,\ s\in [-\ep,+\ep] $. The immersion $\theta_x$ obviously
 extends  as a map (not an immersion) $\Theta_x : B_x\times S^1\to J^1(L,\R)$ 
 by the following formula:
 $$ \Theta_x(z,t)=\theta_x(y,t+s)\,. $$
 For a given $x\in L$, $t\mapsto \Theta_x(-,t)$ is a periodic positive 
contact regular 
 homotopy\footnote{We recall that a regular homotopy is a homotopy through
 immersions.}  of $B_x $ into $J^1(L,\R)$ starting from
 $B_x\hookrightarrow J^1(L,\R)$ at $t=0$ (in fact, if  $B_x$ is small, it is 
an isotopy). 
 The family of the germs of $\Theta_x$ along $\{x\}\times S^1$ may be thought
 of as
a {\it formal solution} of the problem that theorem \ref{theo} solves.
In the next section we follow the book by  Eliashberg-Mishachev \cite{eliash}
for modifying this formal solution into a genuine solution.\\

\section{Towards a genuine solution}

\medskip
In the sequel, it is more convenient to work with a cubication $\mathcal C$ 
(cell decomposition 
made of cubes) instead of a triangulation of $L$. By adding the barycenter of
 each simplex of a triangulation one easily gets a cubication. In the sequel
the germ at $x$ of some periodic homotopy $\Theta:B\times S^1\to J^1(L,\R)$ will
mean the germ of $\Theta$ along $\{x\}\times S^1$.

\begin{prop}\label{1-cell} 
There exist the following families:

\nd {\rm 1)} For each edge $\si\in \mathcal C$ whose end points are $x_0$ 
and $x_1$, there exist a small neighborhood $B_\si$ of $\si$ in $ J^1(L,\R)$
 and  a periodic positive contact regular 
homotopy $\Theta_\si:B_\si\times S^1\to J^1(L,\R)$ starting from $B_\si\hookrightarrow J^1(L,\R)$ at $t=0$ and whose germ at $\{x_i\}$ 
is the one of $\Theta_{x_i}$ for $i=0,1$.

\nd {\rm 2)} There exists a smooth family $\Theta^1_x,\ x\in L,$ of periodic
 positive contact regular homotopies, defined near $x$, whose germ at $x$ is the one of 
$\Theta_\si$ when $x\in \si$.

\end{prop} 

\nd {\bf Proof.} Take the barycentric parametrization of $\si$, 
$\ga(u) $ with  $\ u\in [0,1]$, and a fine subdivision 
$u_0=0,u_1=1/N,\dots,u_N=1$. Let $x_k=\ga(u_k)$, $u'_k=\ds\frac{u_k
+u_{k+1}}{2}$ and $x'_k= \ga(u'_k)$. \\

\nd {\sc Lemma.} {\it  If $N$ is large enough, there exists a periodic 
 Hamiltonian contact isotopy $F^t,\ t\in S^1,$ of
$B_{x_k}$ with support in $2/3B_{x_k}$, being Identity at $t=0$, such that: 
\bi
\item[i)] $\Theta_{x_k}\circ F$ is a positive regular homotopy, where the 
composition is meant at each time $t\in S^1$.
\item[ii)] $x'_k$ belongs to $1/3B_{x_k}$ and the germ of $\Theta_{x_{k+1}}$ 
at $x'_k$ is the one of $\Theta_{x_k}\circ F$.\\
\ei}

\nd {\sc Proof of the lemma}. If $N$ is large enough,
 $\Theta_{x_{k+1}}(x'_k,t)$ belongs to the ball $\Theta_{x_{k}}(1/3B_{x_k},t)$
for every $t\in[0,1]$. Therefore,
 near $\{x'_k\}\times S^1$, $\Theta_{x_{k+1}}$ can be lifted
 to the source of $\Theta_{x_k}$; in other words, the germ of 
 $\Theta_{x_{k+1}}$ along $\{x'_k\}\times S^1$ has a time preserving
 factorization 
through $\Theta_{x_k}$. This lift is a periodic contact Hamiltonian
 isotopy of embeddings of 
a small ball centered at $x'_k$ into $1/3B_{x_k}$. Moreover, if $N$ is large,
 it can be chosen 
$\ep$-close to Identity in the $C^1$-topology. It extends as a Hamiltonian
 isotopy $F$ of $B_{x_k}$ supported in  $2/3B_{x_k}$
and $2\ep$-close to Identity. In general, such an  $F$ is not 
a periodic isotopy; $F^1$ is Identity only on a neighborhood 
$N(x'_1)$ of $x'_1$ and outside 
$2/3B_{x_k}$. But being $C^1$-close to Identity, $F^1$ is isotopic 
to Identity by a contact Hamiltonian
isotopy supported in  $W:=2/3B_{x_k}\setminus N(x'_1) $; indeed, the group 
$Diff_{cont}(W,\partial W)$ is locally contractible\footnote{ 
$Diff(W,\partial W)$ is locally contractible and
 there is a locally trivial fibration
$Diff_{cont}(W,\partial W)\hookrightarrow 
Diff(W,\partial W)\to Cont(W,\partial W)$ whose base is 
the locally contractible space of contact structures coinciding 
with a given one 
along the boundary.}. This allows us 
to modify $F$ so that it becomes a periodic isotopy.
As $\Theta_{x_k}$ is a positive regular
 homotopy, if 
$F$ is close enough to Identity, $\Theta_{x_k}\circ F$ is still 
positive.
${}$ \bull\\

For $0<k<N$, let $B_k:=B_{x_k}$ and let $\Theta'_{k}$ denote the modified
 regular homotopy   described above. We also choose a contact Hamiltonian 
diffeomorphism $\psi_k$ with support in a small neighborhood  of some ray
 $R_k$ in
$B_k$, leaving a neighborhood of $x'_{k-1}$ and $x'_k$  fixed, and moving 
$x_k$ into
$B_k\setminus (2/3B_k)$. We look at the path $\ga'$ defined by:
\bi 
\item[] $\ga'(u)=\psi_k\bigl(\ga(u)\bigr)$ when $0<k<N$ and 
$u\in [u'_{k-1},u'_k]$,
\item[] $\ga'(u)=\ga(u)$ when $u\in [0,u'_0]$ or $u\in [u'_{N-1},1]$ .
\ei
We think of the process changing $\ga$ to $\ga'$ as a {\it making waves } process 
on $\si$, according to the terminology of Bill Thurston in \cite{thurston}.

We now define a positive contact regular homotopy $\Theta'$ 
on a neighborhood of 
$\ga'$ by the following formulas which are matching:
\bi
\item[] $\Theta'\bigl(\ga'(u),t\bigr)=\Theta'_{k}(\ga'(u),t\bigr)$ when
 $0<k<N$ and $u\in [u_k,u'_k]$,
\item[] $\Theta'\bigl(\ga'(u),t\bigr)=\Theta_{k}(\ga'(u),t\bigr)$  
when $0<k\leq N,\ u\in [u'_{k-1},u_k]$ or $k=0,\ u\in[0,u'_0]$.
\ei
This regular homotopy  has the property which is required 
in point 1), except that it is not defined near $\si$ but near 
the path $\ga'$. 

If the rays $R_k$ are chosen to be mutually disjoint and so that
$R_k\cap L=\{x_k\}$ for every $k$, there is  a Hamiltonian contact 
diffeomorphism $\psi$, with compact support in $J^1(L,\R)$, such that 
$\ga'=\psi\circ\ga$ and $\psi$ leaves the other edges of $\mathcal C$ fixed.
So we define for every $t\in S^1$ 
$$\Theta^t_\si=\psi^{-1}\circ\left(\Theta'\right)^t\circ
\psi.$$ It is well-defined on a small neighborhood of $\si$ and 
$t\mapsto \Theta^t_\si$
  is still a positive homotopy as $\psi$ is independent of
 $t\in S^1$. Hence point 1) is proved.\\

\nd 2) When $x\in \si$, it is not difficult to interpolate between 
$\Theta_x$ and the germ of $\Theta_\si$ at $x$ (interpolate between $\psi$ and
 Identity and between
$\Theta_k$ and $\Theta'_k$ if $x=\ga(u),\ u\in [u_k,u'_k]$). 
This interpolation extends 
when $x$ is close to $\si$. Using a partition of unity, 
one easily finds a family 
$\Theta^1_x$ with the desired property.
\bull\\

Applying  proposition \ref{1-cell} above simultaneously to each  1-cell of 
$\mathcal C$ yields a positive contact regular homotopy defined near the
 1-skeleton. So we have  
a periodic contact positive regular homotopy $\Theta^1$ defined near the
 1-skeleton together with a family of $\Theta_x$ defined near each point 
$x\in L$. For going further, as in \cite{eliash}, we need a parametric 
version of \ref{1-cell}.\\

\begin{prop}\label{1-par} Let $\si$ be a 2-cell in $\mathcal C$, $\tau$ be
 a 1-face of $\si$ and $\tau^*$ be a non-parallel 1-face of $\si$; so $\si$ 
 is foliated by  intervals
$\tau(y)$ parallel to  $\tau$, starting at $y\in\tau^*$ and ending at a point
 of the edge opposite to $\tau^*$. Then there exist periodic positive 
contact regular homotopies 
$\Theta_{\tau(y)}$, $y\in \tau^*$, defined near $\tau(y)$ and depending smoothly on $y$. 
Moreover 
its germ at any point of $\partial \si$ is the one of $\Theta^1$.

\end{prop}

\nd {\bf Proof.} Note that $\si $ is a square. Clearly the proof we have done for one edge in \ref{1-cell}
 can be performed with parameters. When $y\in \partial \tau^*$ ($y=0$ or 1),
 as $\Theta^1$ is already defined, it is not necessary to replace $\tau(y)$ 
by a very much oscillating $C^0$-approximation of $\tau(y)$. In other words,
 the contact diffeomorphism $\psi$ from the proof of \ref{1-cell} can be 
interpolated with Identity when $y$ approaches one end point of $\tau^*$.
\bull \\ 
 
\nd {\bf Proof of theorem \ref{theo}.} Here we explain how
 to construct a periodic
positive contact regular homotopy defined near a 2-cell $\si$ of $\mathcal C$
extending $\Theta^1$, the regular homotopy 
we have near the 1-skeleton $L^{[1]}$. When doing it for all 2-cells 
simultaneously,
we get the desired regular homotopy $\Theta^2$ defined near the 2-skeleton 
 $L^{[2]}$.
And one goes on recursively until $\Theta ^n$ which is the desired regular
 homotopy.\\

For constructing $\Theta^2$ near $\si$, we use proposition \ref{1-par} 
which yields a 
1-parameter family $\Theta^1(y)$ of periodic contact regular homotopy
 defined near $\tau(y)$. As in the case of a 1-simplex (proposition 
\ref{1-cell}), we discretize the $y$-interval, deform slightly 
$\Theta^1\left(\tau(y_k)\right)$ so that they glue together and yield a
 periodic 
contact regular homotopy near $\psi(\si)$, where $\psi$ 
is a contact diffeomorphism of 
$J^1(L,\R)$,
{\it making waves} on $\si$. The process changing $\si$ to $\psi(\si)$
 is the 2-dimensional analogue of the one that we have described very precisely 
for a 1-cell; it is a universal process once we know 
  $\Theta^1\left(\tau(y)\right)$ for every $\ y\in\tau^*$. \bull

${}$\\

\end{document}